\documentstyle[fullpage,12pt,amssymb]{amsart}
\input{psfig}
\def\SBIMSMark#1#2#3{
 \font\SBF=cmss10 at 10 true pt
 \font\SBI=cmssi10 at 10 true pt
 \setbox0=\hbox{\SBF Stony Brook IMS Preprint \##1}
 \setbox2=\hbox to \wd0{\hfil \SBI #2}
 \setbox4=\hbox to \wd0{\hfil \SBI #3}
 \setbox6=\hbox to \wd0{\hss
             \vbox{\hsize=\wd0 \parskip=0pt \baselineskip=10 true pt
                   \copy0 \break%
                   \copy2 \break%
                   \copy4 \break}}
 \dimen0=\ht6   \advance\dimen0 by \vsize \advance\dimen0 by 8 true pt
                \advance\dimen0 by -\pagetotal
 \dimen2=\hsize \advance\dimen2 by .25 true in
  \openin2=publishd.tex
  \ifeof2\setbox0=\hbox to 0pt{}
  \else 
     \setbox0=\hbox to 3.1 true in{
                \vbox to \ht6{\hsize=3 true in \parskip=0pt  \noindent  
                \input publishd.tex 
                \vfill}}
  \fi
  \closein2
  \ht0=0pt \dp0=0pt
 \ht6=0pt \dp6=0pt
 \setbox8=\vbox to \dimen0{\vfill \hbox to \dimen2{\copy0 \hss \copy6}}
 \ht8=0pt \dp8=0pt \wd8=0pt
 \copy8
 \message{*** Stony Brook IMS Preprint #1, #2 ***}
}

\newtheorem{thm}{Theorem}[section]
\newtheorem{cor}[thm]{Corollary}
\newtheorem{lem}[thm]{Lemma}
\newtheorem{prop}[thm]{Proposition}

\theoremstyle{definition}

\newcommand{\QED}{\rlap{$\sqcup$}$\sqcap$\smallskip}

\theoremstyle{remark}

\newcommand{\diam}{\operatorname{diam}}

\newcommand{\cl}{\operatorname{cl}}
\newcommand{\tl}{\tilde}
\newcommand{\eps}{\epsilon}
\newcommand{\di}{\partial}
\newcommand{\ra}{\rightarrow}
\def\noi{\noindent}

\newcommand{\id}{\operatorname{id}}
\newcommand{\length}{\operatorname{length}}
\newcommand{\dens}{\operatorname{dens}}

\newcommand{\CC}{\cal C}
\newcommand{\DD}{\cal D}

\newcommand{\II}{\cal I}

\newcommand{\NN}{\cal N}

\newcommand{\RR}{\cal R}

\numberwithin{equation}{section}

\newcommand{\secref}[1]{\S\ref{#1}}
\newcommand{\lemref}[1]{Lemma~\ref{#1}}

\newcommand{\A}{{\Bbb A}}
\newcommand{\C}{{\Bbb C}}
\newcommand{\D}{{\Bbb D}}

\renewcommand{\L}{{\Bbb L}}

\newcommand{\Q}{{\Bbb Q}}
\newcommand{\R}{{\Bbb R}}
\newcommand{\T}{{\Bbb T}}
\newcommand{\V}{{\Bbb V}}
\newcommand{\U}{{\Bbb U}}
\newcommand{\W}{{\Bbb W}}
\newcommand{\X}{{\Bbb X}}
\newcommand{\Z}{{\Bbb Z}}

\newcommand{\h}{{\bf h}}
\newcommand{\la}{{\lambda}}
\renewcommand{\i}{{\bar i}}
\renewcommand{\j}{{\bar j}}
\def\sm{\smallsetminus}
\def\Bf{{\bold{f}}}
\def\Bh{{\bold{h}}}
\def\Bg{{\bold{g}}}
\def\BG{{\bold{G}}}
\def\BT{{\bold{T}}}
\def\Bphi{{\bold{\Phi}}}
\def\Bpsi{{\bold{\Psi}}}
\def\B0{{\bold{0}}}
\renewcommand{\lq}{``}
\renewcommand{\rq}{''}

\renewcommand{\marginpar}[1]{}
\catcode`\@=12

\def\Empty{}
\newcommand\oplabel[1]{
  \def\OpArg{#1} \ifx \OpArg\Empty {} \else
  	\label{#1}
  \fi}
		
\newcommand{\comm}[1]{}

\begin{document}

\hfill {\it You first plow in the dynamical plane and}\hskip 1.5in

\hfill{\it  then harvest in the parameter plane.}\hskip 1in

\hfill{\it  Adrien Douady}
\bigskip\bigskip

\title[parapuzzle]{Dynamics of quadratic polynomials, III\\
                          Parapuzzle and SBR measures.
}
\author {Mikhail Lyubich }
\date{June 25, 1996}
\maketitle
\SBIMSMark{1996/5}{June 1996}{}

\section{Introduction}
This is a continuation of notes on dynamics of quadratic polynomials.
In this part we  transfer the geometric result of \cite{L2}
to the parameter plane. To any parameter value $c\in M$ in the Mandelbrot set
(which lies outside of the main cardioid and little Mandelbrot sets attached to it) 
we associate a \lq{principal  nest of parapuzzle pieces}''
   $D^0(c)\supset D^1(c)\supset\dots$ corresponding 
 to the generalized renormalization type of $c$.  Then we  prove:

\proclaim Theorem A.  The moduli of the parameter annuli  $\mod (D^l(c)\sm D^{l+1}(c))$
grow at least  linearly
(see \S\ref{parapuzzle geometry} for a more precise formulation).

This result was announced at the Colloquium in honor of Adrien Douady 
(July 1995),
and in the survey \cite{L4}, Theorem 4.8.
The main motivation for this work was to  prove the following:

\proclaim Theorem B (joint with Martens and Nowicki). Lebesgue
almost every real quadratic $P_c: z\mapsto z^2+c$ which is non-hyperbolic and at most
finitely  renormalizable has a 
finite  absolutely continuous invariant measure. 

More specifically, Martens and Nowicki [MN]  have given a geometric
criterion for existence of a finite
 absolutely continuous invariant measure (acim) in terms of the ``scaling factors''.
Together with the result of \cite{L5} on the exponential decay of the scaling factors
in the quasi-quadratic case this yields existence of the acim
once  ``the principal nest is eventually free from the central cascades''.
 Theorem A above implies that this condition is satisfied for
almost all real  quadratics  which are
non-hyperbolic and at most finitely  renormalizable (see Theorem \ref{measure zero}).
Note  that Theorem A also implies that this condition is satisfied 
on a  set of positive measure, which yields a new  proof of
Jacobson's  Theorem  \cite{J} (see also Benedicks \& Carleson \cite{BC}).  

A  measure $\mu$ will be  called SBR (Sinai-Bowen-Ruelle) if 
\begin{equation}\label{SBR}
  {1\over n}\sum_{k=0}^{n-1} \delta_{f^k x}\to \mu
\end{equation}
for a set of $x$ of positive {\it Lebesgue } measure. It is known that if an SBR 
measure exists for a real quadratic map $f=P_c$,   $c\in [-2,1/4]$,
on its invariant interval $I_c$,  then it is unique and
(\ref{SBR}) is satisfied for Lebesgue  almost all $x\in I_c$
 (see Introduction of \cite{MN} for a more detailed discussion).
Theorem B yields

\proclaim Corollary. For almost all $c\in [-2,1/4]$, 
  the quadratic polynomial $P_c$
  has a unique SBR measure on its invariant interval $I_c$.

Let us now take a closer look at Theorem A. 
It nicely  fits to the general philosophy
of correspondence between the dynamical and parameter plane. 
This philosophy was introduced to holomorphic dynamics
  by Douady and Hubbard \cite{DH1}. Since then,
there have been many beautiful results in this spirit, see 
  Tan Lei \cite{TL}, Rees \cite{R}, Shishikura \cite{Sh}, Branner-Hubbard \cite{BH},
Yoccoz (see \cite{H}). 

In the last work, special tilings into \lq{parapuzzle pieces}\rq
 of the parameter plane are  introduced. Its  main geometric result is that
the tiles around at most finitely renormalizable points shrink.
   It was done by transferring, in an ingenious way,  the corresponding
dynamical information into the  parameter plane. 

In \cite{L2} we studied the rate at which the dynamical tiles shrink. Our main geometric
result is that   the moduli of the principal nest
of dynamical annuli  grow at least linearly. 
Let us note that the way we transfer this result  to the parameter plane (Theorem~A)
is substantially  different from that of Yoccoz.
 Our main tool is provided by  holomorphic motions
whose  transversal quasi-conformality yields comparison between the dynamical and
parameter moduli.

The properties of holomorphic motions are discussed in  \S\ref{background}.
In \S\ref{parapuzzle combinatorics}
  we describe the principal parameter tilings
according to the generalized renormalization types of the maps. 
In \S\ref{parapuzzle geometry} we prove Theorem~A.
 In the last section, \S\ref{measure}, 
we derive the consequence for the real quadratic family (Theorem~B).

Let us finally draw the reader's attention to the  work of LeRoy Wenstrom \cite{W}
which studies in detail parapuzzle geometry near the Fibonacci parameter value.

\section{Background}\label{background}
\subsection{ Notations and terminology.}
 $\D_r=\{z: |z|<r\}$, $\D\equiv \D_1$,  
$\T_r=\{z: |z|=r\}$, $\A(r,R)=\{r<|z|<R\}$. The closed and semi-closed
annuli are denoted accordingly: $\A[r,R]$, $\A(r,R]$, $\A[r,R]$. 

By a {\it topological disc} we will mean a simply connected domain $D\subset \C$ 
whose boundary is a Jordan curve.

Let $\pi_1$ and $\pi_2$ denote the coordinate projections
$\C^2\rightarrow \C$. Given a set $\X\subset\C^2$, we denote by
$X_\la=\pi_1^{-1}\{\la\}$ its vertical cross-section through $\la$ (the ``fiber" over
$\la$). Vice versa, given a family of sets $X_\la\subset \C$, $\la\in D$, we will use
the notation: $\X=\cup_{\la\in D}X_\la=\{(\la,z)\in \C^2: \la\in D, z\in X_\la\}$.

Let us have a discs fibration  $\pi_1: \U\ra D$ over a topological disc
 $D\subset \C$ (so that the sections $U_\la$ are topological discs, and
the closure of  $U$ in $ D\times\C$ is homeomorphic to $D\times\bar\D$ over $D$).  
In this situation we call $\U$  an (open) {\it topological bidisc} over  $D$. 
We say that this fibration admits an extension to the boundary $\di D$
if the closure $\bar U$  of $\U$ in $\C^2$ is homeomorphic
over $\bar D$  to $\bar D\times \bar \D$ .
 The set $\bar \U$ is called a (closed) bidisc.
 We keep the notation $\U$ for the fibration of {\it open}
discs over the closed disc $\bar D$ (it will be clear from the context over which
set the fibration is considered). 

If $U_\la\ni 0$, $\la\in D$, we denote by $\bold{0}$ the zero section of the fibration. 

If the fibration  $\pi_1$ admits an extension over the boundary $\di D$, 
we define  the {\it frame} $\delta \U$ as the topological torus 
$\cup_{\la\in \di D} \di U_\la$. A section $\Bphi: D\ra \U$ is called {\it proper}
if it is continuous up to the boundary and $\Bphi(\di D)\subset\delta U$.

We assume that the reader is familiar with the theory of quasi-conformal maps
(see e.g., \cite{A}). 
We will use a common abbreviation $K$-qc for \lq{$K$-quasi-conformal}\rq.  
Notation $a_n\asymp b_n$ means, as usual, the the ratio $a_n/b_n$ is 
positive, and bounded away from 0 and $\infty$.

\subsection{Holomorphic motions} Given a domain $D\subset \C$ with a base point $*$
 and a set  $X_*\subset\C$, a {\it holomorphic motion} $\Bh$ of $X_*$ over $D$ is a family
of injections $h_\la: X_*\ra\C$, $\la\in D$, such that $h_*=\id$ and $h_\la(z)$ is
holomorphic in $\la$ for any $z\in X_*$.  We denote $X_\la=h_\la X_*$.   Let us summarize
fundamental properties of holomorphic motions which are usually referred to as the
$\la$-{\it lemma}. It consists of two  parts: extension of the motion and transversal
quasi-conformality, which will be stated separately.

The consecutively  improving versions of the Extension Lemma appeared in
 [L1] and [MSS], [ST],[BR],[Sl]. 
 The final result, which will be actually exploited below, is due to Slodkowsky:

\proclaim Extension Lemma.  A holomorphic motion $h_\la: X_*\ra X_\la$  of a set
$X_*\subset\C$ over a  topological disc $D$ admits an extension   to a holomorphic motion
$H_\la: \C\ra\C$ of the whole complex plane over $D$.

\proclaim Quasi-Conformality Lemma [MSS].
Let $h_\la: U_*\ra U_\la$ be a holomorphic motion of a domain $U_*\subset\C$
over a hyperbolic domain $D\subset\C$. Then the maps $h_\la$ are
 $K(r)$-quasi-conformal, where $r$ is the hyperbolic distance between $*$ and $\la$ in $D$.

A holomorphic motion  $h_\la: U_*\ra U_\la$ over $D$
 can be viewed as a complex one-dimensional foliation of the domain
 $\U=\cup_{\la\in D} U_\la$,
whose leaves are graphs of the functions  $\la\mapsto h_\la(z)$, $z\in U_*$.  A {\it
transversal} to the motion is a complex one dimensional submanifold of $\C^2$ which
transversally  intersects every leaf at one point
(so that \lq{transversal}" means {\it global} transversal).  
Given two transversals $X$ and $Y$, we thus have a well-defined holonomy map
$H: X\ra Y$, $H(p)=q$ iff $p$ and $q$ belong to the same leaf.

\begin{cor}\label{transversal qc structure}  Holomorphic motions are locally transversally
quasi-conformal. More precisely,
  for any two transversals $X$ and $Y$, the holonomy map $H: X\ra Y$
  is locally  quasi-conformal. If $H(p)=q$ then the
   local dilatation of $H$ near $p$ depends only on the hyperbolic distance between the
 $\pi_1(p)$ and $\pi_1(q)$ in $D$.
\end{cor}

\begin{pf}  Let $p=(\la,\alpha)$, $q=(\mu,\beta)$. By the $\la$-Lemma, the map
$G=h_\mu\circ h_\la^{-1}: U_\la\ra U_\mu$ is quasi-conformal, with dilatation depending only
on the hyperbolic distance between $\la$ and $\mu$ in $D$. 
 Hence a little  disc $D(\alpha,\eps)\subset U_\la$ is mapped by $G$
 onto an  ellipse  $Q_\eps\subset U_\mu$
with bounded eccentricity about $\beta$ (where the bound depends only on the hyperbolic
distance between $\alpha$ and $\beta$).

 But the holonomy
 $U_\la\ra X$ is asymptotically conformal near $p$. To see this,
 let us select a  holomorphic coordinates $(\theta,z)$ near $p$ 
in such a way that $p=0$ and the leaf via $p$ becomes the parameter
axis. Let $z=\psi(\theta)=\eps+\dots$ parametrizes  a nearby
 leaf of the foliation,
while $\theta=g(z)=bz+\dots$ parametrizes the transversal  $X$.
 
 Let us do the rescaling $z=\eps\zeta, \theta=\eps\nu$. In these new coordinates, the
above leaf is parametrized by the function
$\Psi(\nu)=\eps^{-1} \psi(\eps\nu)$,
$|\nu|<R$, where $R$ is any  fixed parameter. 
Then $\Psi'(\nu)=\psi'(\eps\nu)$ and $\Psi''(\nu)=\eps
\psi''(\eps\nu)$.  Since the family of functions $\{\psi(\nu)\}$ is normal,
 $\Psi''(\nu)=O(\eps)$. 
  Moreover,  $\psi$ uniformly goes to 0 as $\psi(0)\to 0$. 
Hence $\Psi'(0)=\psi'(0)\leq \delta_0(\eps)$,
where $\delta_0(\eps)\to 0$ as $\eps\to 0$.  Thus
$\Psi'(\nu)=\delta_0(\eps)+O(\eps)\leq\delta(\eps)\to 0$ as $\eps\to 0$ uniformly  for all
$|\nu|<R$.
It follows that  $\Psi(\nu)=1+O(\delta(\eps))=1+o(1)$ as $\eps\to 0$. 

On the other hand, the
the manifold $X$ is parametrized in the rescaled coordinates
 by a function $\nu=b\zeta+0(1)$. Since the transversal intersection persists,
 $X$ intersects  the leaf at the point $(\nu_0,\zeta_0)=(1,b)(1+o(1))$
(so that $R$ should be selected bigger than $b$).
In the old coordinates the intersection point is $(\theta_0,z_0)=(\eps, b\eps)(1+o(1))$.

 Thus the holonomy from $U_\la$ to $X$ transforms the disc of radius
$|\eps|$ to an ellipse with small eccentricity, which means that this holonomy is
asymptotically conformal. As the holonomy from $U(\mu)$ to $Y$ is also asymptotically
conformal, the conclusion follows.   
\end{pf}

\smallskip \noindent{\bf Remark.}
  The above argument also  shows that if the motion is $K$-qc over the
whole domain $D$ (i.e, all maps $h_\la$ are $K$-qc, $\la\in D$)
then it is transversally $K$-qc.

\subsection{Winding number}\label{winding} Given two curves $\psi_1, \psi_2: \di D\ra \C$ 
 such that $\psi_1(\la)\not=\psi_2(\la),\;\la\in\di D$,
 we can define the winding number
 of the former about the latter  as the increment of 
${1\over 2\pi} \arg(\psi_1(\la)-\psi_2(\la))$ as $\la$ wraps once around $\di D$. 

Let us have a bidisc  $ \U$ over $\bar D$. 
Given a proper section $\Bphi: D\ra \U$ let us define its {\it winding number} as
follows.
Let us  mark the torus
$\delta\U$ with the homology basis $\{[\di D], [\di\U_*]\}$.
 Then the winding number
$w(\Bphi)$ is the second coordinate of the curve 
$\Bphi:\di D\rightarrow \delta\U$ with respect to this basis.

\proclaim Argument Principle. 
 Let us have a bidisc $\U$ over $\bar D$ and a proper
holomorphic section  $\Bphi: D\ra \U$, $\phi=\pi_2\circ\Bphi$.  Let
$\Bpsi: \bar D\ra\U$ be  another continuous section 
 holomorphic in $D$, $\psi=\pi_2\circ\Bpsi$.
 Then the number of solutions of the equation
$\phi(\la)=\psi(\la)$ counted with multiplicity is equal to the winding number $w(\Bphi)$.

\noindent{\bf Proof.} Indeed, $w(\Phi)$ is equal to the winding number of
$\phi$ around $\psi$, which is equal, \newpage
by the standard Argument Principle, to the
number of roots of the equation  $$\phi(\la)=\psi(\la).$$
\comm{
 Let us consider any homotopy 
$$\psi_t: \di D\ra \C,\quad 0\leq t\leq 1,$$   between $\psi_0(\la)\equiv 0$ and $\psi_1=\psi$,
such that $\psi_t(\la)\in U_\la$. Let us define the winding number
$w_t$ of $\phi$ around $\psi_t$ as the increment of 
$(1/2\pi)\arg(\phi(\la)-\psi_t(\la))$ as $\la$ wraps once around
$\di D$.  Since the curves $\la\mapsto \phi(\la)-\psi_t(\la)$ don't pass through the origin, the
winding number $w_t$ stays constant. Thus $w_0=w_1$. But by the standard
Argument Principle
 $w_0=w(\Bphi)$, while  $w_1$ 
is equal to the number of solutions of the equation
$\phi(\la)=\psi(\la)$. \QED  end comm}

\section{Parapuzzle combinatorics}\label{parapuzzle combinatorics}

\subsection{Holomorphic families of generalized quadratic-like maps}\label{quadratic-like families}

 Let us consider a topological disc $D\subset\C$ with a base point $*\in D$,
and a family of topological bidiscs
$\V_i\subset \U\subset\C^2$ over $D$ ({\it tubes}),
such that  the $\V_i$ are pairwise disjoint.
\comm{Let $V_{i,\la}=\pi^{-1}(\la)\cap \V_{i,\la}$ stand for the vertical
cross-sections of the $\V_{i}$ through $\la\in D$
(the fibers over $\la$), and $U_\la$ have the similar meaning for $\U$.}
 We assume that $V_{0,\la}\ni 0$.

Let
 \begin{equation}\label{g}
\Bg: \cup\V_i\rightarrow \U
\end{equation}
 be a fiberwise
map, which admits a holomorphic extension to some neighborhoods of the
$\V_i$ (warning: these extensions don't fit), and
whose fiber restrictions
$$\Bg(\la,\cdot)\equiv g_\la: \bigcup_i V_{i,\la}\rightarrow U_\la,\quad \la\in D,$$
are generalized quadratic-like maps with the critical point at 
$0\in V_\la\equiv V_{0,\la}$ (see \cite{L2}, \S 3.7 for the definition).
We will also assume that the discs  $U_\la$ and $V_{i,\la}$ are bounded by quasi-circles.

Let us also assume that there is a holomorphic motion $\Bh$, 
\begin{equation}\label{h}
h_\la: (\di U_{*}, \bigcup_i\di V_{i,*})\rightarrow
           (\di U_{\la}, \bigcup_i\di V_{i,\la}),
\end{equation}
with a base point at $*\in D$,  which respects action of $\Bg$:

\begin{equation}\label{respect dynamics}
h_\la\circ g_{*}(z)=g_\la\circ h_\la(z), \quad for \quad z\in
\cup \di V_{i,*}.
\end{equation}
A {\it holomorphic family $(\Bg,\Bh)$ of (generalized) quadratic-like maps over $D$} is
a map (\ref{g}) together with a holomorphic motion (\ref{h}) satisfying 
(\ref{respect dynamics}). We will sometimes reduce  the notation to $\Bg$.
In case when the domain of $\Bg$ consists of only one tube $\V_0$, we refer to $\Bg$
as {\it DH quadratic-like family} (for \lq{Douady and Hubbard}\rq, compare \cite{DH2}).

\bigskip
\centerline{\psfig{figure=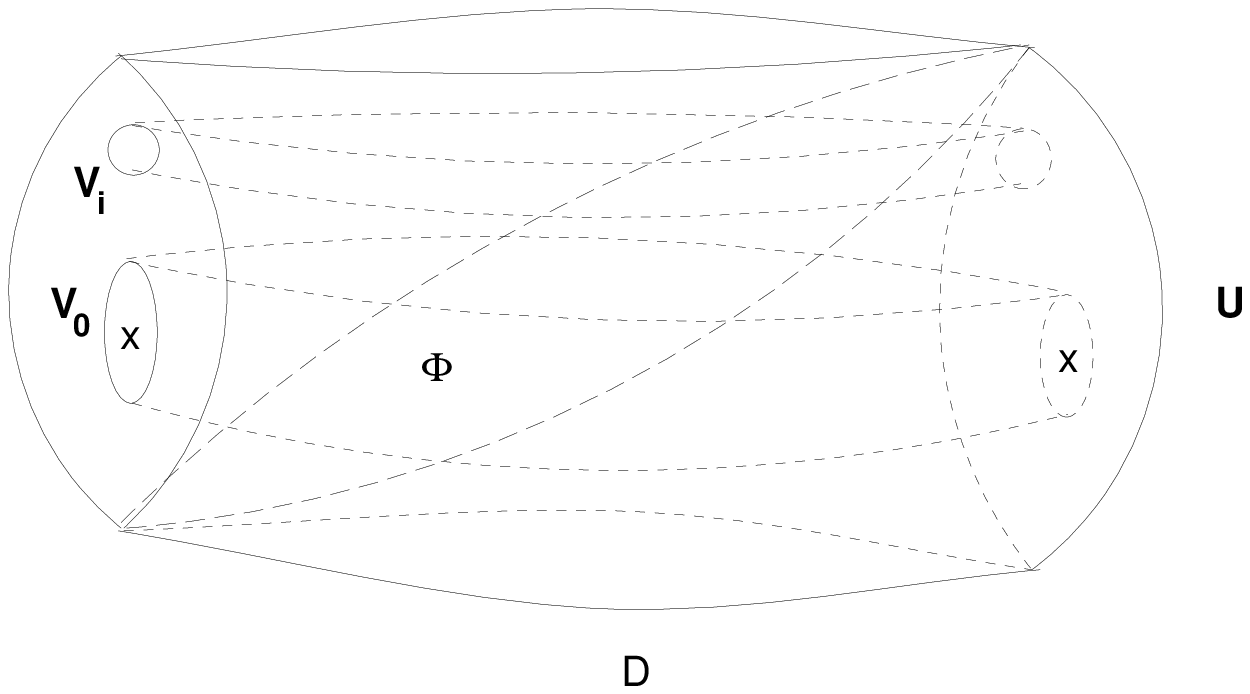,width=.8\hsize}}
\bigskip\centerline{Figure 1: Generalized quadratic-like family.}
\bigskip

Note that by the Extension Lemma we can always assume that the holomorphic motion
(\ref{h}) is extended to the whole disc $U_*$ (though without respect of dynamics
outside $\cup\di V_{i,*}$). 

Let us now consider the critical value function $\phi(\la)=g_\la(0),$
$\Bphi(\la)=\Bg(\la,0)\equiv (\la,\phi(\la))$.
Let us say that $\Bg$ is a {\it proper} (or {\it full}) holomorphic family if
the fibration $\pi_1: \U\ra D$ admits an extension to the boundary  $ \bar D$,
$\bar V_i\subset U$, and $\Bphi:  D\ra \U$ is a proper section.
 Note that  the fibration $\pi_1: V_0\ra D$ cannot be extended to $\bar D$,
 as the domains $V_{\la,0}$ pinch to figure eights as $\la\to \di D$. 

Given a proper holomorphic family $\Bg$ of generalized quadratic-like maps, 
let us define its {\it
winding number} $w(\Bg)$ as the winding number of the critical value $ \phi(\la)$ about
the critical point 0.  By the Argument Principle, it is equal to the winding number
of the critical value about any section  $\bar D\ra \U$.

We will also  face the situation when $\Bg$ does not map every tube $\V_{i}$
onto the whole tube  $\U$ but rather on some other tube $\V_{j}$,
 while all the rest properties listed above are still valid
(see \S\ref{initial family}).
Then we call $\Bg$ a holomorphic family of Markov maps. 

Let $\mod(\Bg)$ stand for the $\inf_{\la\in D}\mod (U_\la\sm V_{0,\la})$.

\subsection{Douady \& Hubbard quadratic-like families }\label{DH families}
 Let us have a proper holomorphic family $\Bf:\V\ra\U$ of DH quadratic-like maps, 
 with winding number 1.
The Mandelbrot set $M(\Bf)$ is defined as the set of $\la\in D$ such that
the Julia set $J(f_\la)$ is connected.
 
Since the $U_\la$ and $V_{i,\la}$ are bounded by quasi-circles,
there is a qc straightening  
 $\omega_*: \cl(U_*\sm V_*)\ra\A[2,4]$ conjugating $f: \di V\ra \di U$  to $z\mapsto z^2$ on $\T_2$.
 The  holomorphic motion $\Bh$ 
 (extended to the whole \lq{condensator}\rq  $\U\sm \V$)
 spreads  this straightening over the whole parameter region $D$. 
We obtain a family of quasi-conformal homeomorphisms
\begin{equation}\label{omega}
 \omega_\la: \cl(U_\la\sm V_\la)\ra \A[2,4] 
\end{equation}
conjugating $f_\la|U_\la$ to $z\mapsto z^2$ on $\T_2$.
 Pulling them back, we obtain for every $f_\la$ the straightening
 $\omega_\la: V_\la\sm \Omega_\la\ra \A(\rho_\la,4)$ 
 well-defined up to the critical point level $\rho_\la=|\omega_\la(0)|$
 (so that for $\la\in M(\Bf)$ it is well-defined on the whole complement of the Julia set).

\smallskip\noindent{\bf Remark.} The motion $\Bh$ can usually be selected
uniformly qc (actually smooth) on $\U\sm \V$ over $D$.  Then the
straightenings $\omega_\la$ are also uniformly qc. 
\smallskip

Let us now define  a map $\xi: D\sm M(\Bf)\ra \A(1,4)$ in the following way:
\begin{equation}\label{parametrization}
\xi(\la)=\omega_\la(f_\la 0).
\end{equation}

\begin{lem}\label{parameter coordinates}  Let $\Bf$ be a DH quadratic-like family with
winding number 1. Then formula \ref{parametrization}  determines a
homeomorphism  $\xi: D\sm M(\Bf)\ra \A(1,4)$. If the holomorphic motion
$\Bh$ is selected  uniformly  qc on the annulus tube $\U\sm \V$ then
$\psi$  is quasi-conformal.
\end{lem}

\begin{pf} Without loss of generality we can assume that $*\in M(\Bf)$. 
Let  $V^n_*=f_*^{-n} U_*,\; n=0,1,\dots$. 

 Let us consider the critical value  graph 
 $C=\Phi(\la)\equiv \{(\la, f_\la 0),\; \la\in D\}$.
 By the Argument Principle, it
intersects at a singe point each leaf of the holomorphic motion $\Bh$ on $\U\sm \V$, so that
the holonomy $\gamma: U_*\sm V_*\ra X$ is a homeomorphism onto the image $R_1$. 
Hence $A_1\equiv \pi_1 R_1\subset D$ is a topological annulus, and the map
$$\xi^{-1}=\pi_1\circ\gamma\circ\omega_*^{-1}: \A[2,4)\ra A_1$$ is a homeomorphism.

Let $\Gamma_1$ be the inner boundary of $A_1$, 
and $D_1$ be the topological disc bounded
by $\Gamma_1$. Since the critical value $f_\la(0)$, $\la\in D_1$, 
does not land
at the leaves of  holomorphic motion $\Bh|D_1$, 
it can be lifted by $\Bf$ to a holomorphic
motion $\Bh_1$ of the annulus $V_*^1\sm V_*^2$ over $D_1$. 
Since the graph $C$ intersects
every leaf belonging to $\di\V^1$ at a single point, the family
 $(\Bf: \V^2\ra\V^1,\Bh)$ is proper over $D_1$ and  has winding number~1. 
Let $A_2=\Phi^{-1}(\V^1\sm \V^2)$.
Then  the same argument as above shows that the map $\xi^{-1}: \A[\sqrt{2},2]\ra A_2$
is also a homeomorphism.

Continuing in the same way, 
we will inductively  construct a sequence of holomorphic motions $\Bh_n$
over nested discs $D_n$, and a nest of adjoint  annuli
  $A_n=D_{n-1}\sm D_n$ which are homeomorphically  mapped by 
$\xi$ onto the round annuli
 $\A[2^{1/(n-1)}, 2^{1/(n-2)}]$. Altogether this shows that $\xi$ is a homeomorphism.

Finally, assume  $\Bh$ is  $K$-qc. Since all further motions
$\Bh_n$ are holomorphic lifts of $\Bh$ over $D_n$ by $\Bf^n$,
 they are $K$-qc over their
domains of definition as well. By Corollary \ref{transversal qc structure} (and the Remark 
afterwards), they are  transversally $K$-qc. 
Moreover, the straightening $\omega_*: U_*\sm J(f_*)\ra \A(1,4)$
 is qc, while the projection $\pi_1: X\ra D$ is conformal. 
Since  $\xi$  is the composition of the straightening,
the  holonomy and the projection,
it is quasi-conformal.
\end{pf} 

\noi{\bf Example} (see \cite{DH1}).
 Let us consider the Mandelbrot set $M$ of
the quadratic family  $P_c: z\mapsto z^2+c$. Let 
$R: \C\sm M\ra \C\sm \D$ be the Riemann mapping. Recall that 
parameter equipotentials and external rays are defined as 
the $R$-preimages of the round circles and radial rays.
Let $\Omega_r$ be the topological disc bounded by the equipotential
 $R^{-1}\{re^{i\theta}: 0\le \theta \le 2\pi\}$
of level $r>1$.

For every $c\in \Omega_4$, let us consider the quadratic-like map
$P_c: V_c\ra U_c$ where $V_c$ and $U_c$ are topological discs bounded by the 
dynamical equipotentials of level 2 and 4 correspondingly.
Then the conformal map $\omega_C: U_c\sm V_c\ra \A(2,4)$ 
conjugates $P_c|\di V_c$ to $z\mapsto z^2 $ on $\T_2$, so that it can 
serve as a straightening (\ref{omega}). With this choice of straightening,
the parameter map $\xi: D\sm M\ra \A(1,4)$ constructed in 
\lemref{parameter coordinates} just coincides with the Riemann map  $R$.
\QED \smallskip

With \lemref{parameter coordinates}, we
 can  extend the  the notion of parameter rays and equipotentials 
to quadratic-like families as the $\xi$-preimages of the polar
coordinate curves in $\A(1,4)$.  
If $\xi(\la)=re^{i\theta}$ then $r$ and $\theta$ are called
the {\it equipotential level} and the {\it external angle} of the parameter value $\la$.
Note that $\di D$ becomes the equipotential of level 4.

\subsection{Wakes and initial Markov families}\label{initial family}
Recall that every quadratic-like map $f: V\ra U$ is hybrid equivalent to a
quadratic polynomial $P_c: z\mapsto z^2+c$ 
(The Straightening Theorem \cite{DH2}). It is constructed
by gluing $f|U$ to $z\mapsto z^2$ on $\C\sm \D_2$,
and pulling  the standard conformal structure on $\C\sm\D$ back to
$U\sm K(f)$ by iterates of $f$. The construction 
depends on the choice of
a qc straightening $\omega: \cl(U\sm V)\ra \A[2,4]$
conjugating $f|\di V$ to $z\mapsto z^2$ on $\T_2$. However, if the
Julia set $J(f)$ is connected, the parameter value $c\equiv \chi(f)$
is  determined uniquely.

Given a quadratic-like family $f_\la: V_\la\ra U_\la$ over $D$
with winding  number 1, let  us consider a family of straightenings
(\ref{omega}) and the corresponding family of quadratic polynomials
$P_{\chi(\la)}: z\mapsto z^2+\chi(\la)$. The following statement
follows from the main result of \cite{DH2}:

\begin{lem}\label{continuity}
Under the circumstances just described,
the straightening $\chi: (D, M(\Bf))\ra (\Omega_4, M)$ is
a homeomorphism of the disc $D$ onto a neighborhood $\Omega_4$
of the Mandelbrot set $M$ bounded by the parameter equipotential of
level 4.
\end{lem}

Note that the above explicit description of the image of $\chi$
follows from the following formula:
\begin{equation}\label{relation}
\xi=R\circ\chi,
\end{equation}
where $\xi $ is defined in (\ref{parametrization}), and $R$ is the Riemann
mapping on the complement of the Mandelbrot set. This formula, in turn,
follows from the definitions of $\xi$ and $\chi$ and the description of $R$
given (see the Example in the end of \S \ref{DH families}).

\lemref{continuity} shows  that the landing properties of the parameter
rays in a quadratic-like family coincide with the corresponding properties
in the quadratic family. This allows us 
 to extend the notions of the parabolic and  Misiurewicz wakes from the quadratic
to the quadratic-like case. Namely, the $q/p$-{\it parabolic wake} $P_{q/p}=P_{q/p}(\Bf)$
is the parameter region in
$D$ bounded by the external rays landing
 at the $q/p$-bifurcation point $b_{q/p}$ on the main cardioid of $M(\Bf)$ and the appropriate
arc of $\di D$.
Dynamically it is specified by the property that for $\la$ in this wake there are 
$p$ rays landing at the $\alpha$-fixed point $\alpha_\la$ of
$J(f_\la)$, and they form a cycle with rotation number $q/p$.

The maps
\begin{equation}\label{satillite family}
f_\la^p: V_\la\ra  U_\la
\end{equation}
   restricted to appropriate domains
form  a quadratic-like family  over the wake (see \cite{D}, \cite{L2}, \S\S 2.5, 3.2).
(The domain $V_\la$ is a thickening of the puzzle piece $Y^{(1+p)}_\la$
 bounded by two pairs of rays landing at the
$\alpha$-fixed and co-fixed points and two equipotential arcs. The domain $U_\la$ is
a thickening of the first puzzle piece $Y^{(0)}_\la$ bounded
by two rays landing at the $\alpha$-fixed point and an equipotential arc.) 
Note however that this family fails to be
proper as the domains $U_\la$ don't admit continuous extension at the root.  

\begin{prop}[see \cite{D}]\label{initial rotation number}
 Let $\Bf$ be a DH quadratic-like family
with winding number 1. Then the winding number of the critical value
$\la\mapsto f_\la^p(0)$ about 0 when $\la$ wraps once about the boundary of the
parabolic wake $\di P_{q/p}$ is also equal to 1.
\end{prop}

By \cite{DH2,D}, the quadratic-like family \ref{satillite family} generates
a homeomorphic copy  $M_{q/p}=M_{q/p}(\Bf)$ of the
 Mandelbrot set attached to the bifurcation point $b_{q/p}$.
 Its complement  $M\sm M_{q/p}$ consists of 
a component containing the main cardioid
and infinitely many {\it decorations} 
(using terminology of Dierk Schleicher) $D_{q/p}^{\sigma,i}$, where $\sigma$ is
a dyadic sequence of length $|\sigma|=n-1$,
 $n=1,2,\dots$,  $i=1,\dots,p-1$. A decoration $D_{q/p}^{\sigma,i}$ touches $M_{q/p}$ at a
Misiurewicz point $\mu=\mu_{q/p}^\sigma$ for which 
$$ f_\mu^{pk}(0)\in Y_\mu^{(1+p)},\; k=0,\dots,n-1,\;\; \text{while}\;\;
      f_\mu^{pn} (0)=\alpha_\mu',$$
where $\alpha_\mu'$ is the $\alpha$-co-fixed point 
(i.e., the $f_\mu$-preimage of the fixed point $\alpha_\mu$). (Such Misiurewicz points
are  naturally labelled by the dyadic sequences).  

Every decoration $D_{q/p}^{\sigma,i}$ belongs to the {\it Misiurewicz wake}
 $\tl O_{q/p}^{\sigma,i}$ bounded by
two parameter rays landing at $\mu_{q/p}^\sigma$ 
(there are $p$  rays landing at this point). 
Let us truncate such a wake by the equipotential of level $4/(pn-1)$. 
We will obtain the
initial puzzle pieces $O^{q/p}_{\sigma,i}$ 
which sometimes will also be called \lq{wakes}\rq. 
They can be  dynamically specified in terms of the initial puzzle 
 (see \cite{L3},\S 3.2). 
Namely, there are $p-1$ puzzle pieces $Z^{(1)}_i$, $n=1,\dots, p-1$,
attached to the co-fixed point $\alpha'$. Pulling them back by $(n-1)$-st iterate
of the double covering  $f^p: Y^{(1+p)}\ra Y^{(1)}$, we obtain $2^{n-1}$ puzzle pieces
$Z^{1+(n-1)p}_{\sigma,i}$ labelled by the dyadic sequences.
The wake $O_{q/p}^{\sigma,i}$ is specified by the property
that $f_\la^{p}0\in Z^{(1+(n-1)p)}_{\sigma,i}$.

By {\it tiling} we will mean a family of topological discs with
disjoint interiors. 
Let us consider the initial tiling constructed in \cite{L3}, \S 3.2:
\begin{equation}\label{initial tiling} 
Y^{(0)}_\la\supset V^0_\la\cup X_{i,\la}^k\cup Z_{j,\la}^{(1+kp)}.
\end{equation}
Let us recall that
$\Z_j^{(1)}$ means $\cup_\la Z_{j,\la}^{(1)}.$

\begin{lem}\label{initial partition} 
Let $\Bf$ be a DH quadratic-like family with winding number 1.
The initial tiling (\ref{initial tiling}) moves  holomorphically
 within the Misiurewicz wake $O=O_{p/q}^{\sigma,i}$. 
The critical value of the return map,
 $\Bphi: O\ra \Z_j^{(1)},\; \Bphi(\la)=  f_\la^{pn}0$, 
is a proper map with winding number 1 (where $n=|\sigma|+1)$.  
\end{lem}

\begin{pf} Indeed, all puzzle pieces of this initial tiling 
are the pullbacks of 
$Z^{(1)}_{j,\la}$. As $\la$ ranges over the wake $O\equiv O_{q/p}^{n,i}$, 
the corresponding  iterates of 0 don't cross the boundary of $Z^{(1)}_{j,\la}$. 
It follows that the boundary of the initial tiling moves holomorphically.

Moreover,  the torus $\delta\Z^{(1)}_j$ 
is foliated by the curves with the same
external coordinates, and one curve corresponding 
to the motion of the $\alpha$-co-fixed point.
By the definition of the Misiurewicz wake, the critical value
$\Bphi(\la)$  intersects   once  every leaf of this foliation when $\la$
  wraps once around $\di O$. Hence
$\Bphi: O\ra\Z_j^{(1)}$ is a proper map with winding number 1.  
\end{pf}

Thus we have the initial Markov partition moving  holomorphically
  over the corresponding Misiurewicz wake. 
The wake $O_{q/p}^{\sigma,i}$ containing a point $\la$ will also be denoted 
by $O(\la)$. 

\subsection{First generalized quadratic-like family}
\label{first generalized family} 
Let us consider a proper DH 
quadratic-like family $\Bf=\{f_\la\}$  over $D$ with winding number 1.
Fix a Misiurewicz wake $O$ of this family. 
 The first generalized
quadratic-like map $g_{1,\la}: \cup V^1_{i,\la}\ra V^0_\la$  is defined as the first return map
to $V^0_\la$ (see \cite{L3}, \S 3.5).  The itinerary of the critical point via the elements of
the initial  tiling (\ref{initial tiling})
determines the parameter tiling $\DD^0$ of a  Misiurewicz wake $O$  by
the corresponding puzzle pieces. Let $D_0(\la)$ stand for such parapuzzle piece containing $\la$.

More precisely,  for any $\la\in O$, let us consider the first landing map to 
$T_\la: \cup L_{\bar i,\la}\ra V_{0,\la}$ (see \cite{L3}, \S 4.4). 
The puzzle piece $L_{\bar i,\la}$ 
is specified by its itinerary $\bar i=(i_0,\dots,i_{l-1})$ 
 under iterates of $f_\la^p$ through non-central
pieces of the initial tiling until the first landing at $V_{0,\la}$:
$$L_{\i,\la}=\{z: f_\la^{pk}z\in V_{i_k},\; k=0,\dots, l-1,\;\; f_\la^{pl}z\in V_0\}.$$
 Moreover, $T_\la$
univalently maps $L_{0,\la}$ onto $V_{0,\la}$.  Let $\i_\la$ stand for the itinerary
of the critical value $\psi(\la)=f_\la^p 0$ through the initial tiling,
 so that $f_\la^p 0\in L_{\i_\la,\la}\equiv Q_\la$.  Then the parapuzzle pieces of the tiling
$\DD^0$ are defined as follows:
 $$D_0(\la)=\{\mu\in O: Q_\mu = Q_{\la}\}.$$

Since the critical value $\psi(\la)$ does not intersect the boundary of the initial tiling
as $\la$ ranges over the wake $O$, the pieces $L_{\bar i, \la}$ form the 
tubes $\L_{\bar i}$ with holomorphically moving boundary. Since the winding number
of  $\Bpsi(\la)=(\la,\psi(\la))$ 
 about the tubes of the initial tiling (\ref{initial tiling}) over $O$ is equal to 1
(by \lemref{initial partition}),
the function $\Bpsi: D_0(\la)\ra \Q$ is proper with
 winding number  1.
Since the first landing map is a fiberwise diffeomorphism of  every tube 
$\L_{\bar i}$ onto $\V_0$, the function $\la\mapsto g_\la(0)=T_\la\circ \psi(\la)$,
$D_0(\la)\ra \V_0$, is also proper with
winding number  1. Thus we have:

\begin{lem}\label{first family}
 Let $\Bf$ be a DH quadratic-like family with 
winding number 1.  
 Then the  first generalized renormalization $\Bg=\{g_\la\}$ is a proper
family with winding number 1.
\end{lem}

 \subsection{Renormalization of holomorphic families.}\label{renormalization of families}
Let us now  have a holomorphic quadratic-like family,
 $\Bg:\cup\V_i\ra\U $ over $D$, see \secref{quadratic-like families}.
Let $\II$ stand for the labeling set of tubes $\V_i$. Remember that
 $\II\ni 0$ and $\V_0\ni\B0$. Let
${\II_\#}$ stand for the set of all finite sequences
$\i=(i_0,\dots, i_{l-1})$ of non-zero symbols
$i_k\in \II\setminus\{0\}$. For any $\i\in \II_\#$, there is a tube
$\V_\i$ such that
$$\Bg^k \V_\i\subset \U_{i_k},\; k=0,\dots, l-1,\quad and \quad 
g^l \V_\i=\U.$$ We call $l=|\i|$ the rank of this tube.
The  map $\Bg^l: V_\i\ra\U$ is a holomorphic diffeomorphism which
fibers over $\id$, that is,
$g^l_\la V_{\i,\la}=U_\la,\; \la\in D$.

Let us now pull the  holomorphic motion $\h$ back to the boundaries
$\di \V_\i$:
$$g_\la^l\circ  h_{\i,\la}(z)=h_\la (g^lz),
\; z\in \di V_{\i,*}.$$
This holomorphic motion is an extension of $\Bh$, as by (\ref{respect
dynamics}) it coincides with $\Bh$ on $\di \U_{i}$. In what follows we
will keep notation $\Bh$ for this extended motion.

Let $\L_{\i}\subset \V_{\i}$ be such a tube that
$\Bg^{l} \L_{\i}=\V_{0}$ where $l=l_\i=|\i|$. Extend the holomorphic motion 
$h_\la$ to the boundaries of these tubes by pulling it back
from $\di\V_0$ to the $\di\L_{\i}$ by $\Bg^{l}$. 

 Let
$\i_\la$  be the itinerary of the critical value $\phi(\la)=g_\la 0$ 
under iterates of $g_\la$ through
the  domains $V_{i,\la}$, until its first return to $V_{0,\la}$.
In other words, let $g_\la(0)\in V_{\i_\la,\la}$ and $g^{l_\la}(0)\in V_{0,\la}$,
where
$l_\la=|i_\la|$. Let $Q_\la=L_{i_\la,\la}$ be the corresponding puzzle piece. 
 
Let us now define a new parameter domain $D'$ as the component of the set
$$\{\la: \phi(\la)\in V_{i_*,\la}\}=\{\la: Q_\la=Q_*\}$$ containing $*$. For $\la\in D'$,
the itinerary of the critical value under iterates of $g_\la$ until the first
return back to $V_{0,\la}$ is the same as for $g_*$ (that is, $\i_*$).
Let us define
 new tubes $\V_j'\subset \V_0$ as the components  of
$(\Bg|\V_0)^{-1} (\L_\i\cap \pi^{-1} D')$.
Let 
\begin{equation}\label{renorm}
\Bg':\cup \V_j'\rightarrow \V_0
\end{equation}
 be the first return map of the union of these tubes onto $\V_0$.

For $\la\in D'$, the critical value $\Bphi(\la)$ does not intersect the boundaries of the
the tubes $\L_\i$. Hence we can pull back the holomorphic motion of $\di\L_\i$ to a 
holomorphic motion $\Bh'$ over $D'$ of $\di\V_j'$.  
Then $(\Bg',\Bh')$ is a holomorphic family over $D'$ which will be called the 
{\it generalized  renormalization} of the family $(\Bg,\Bh)$.

If $\Bg$ is a proper family then $\Bg'$ is clearly proper as well.
Moreover, $w(\Bg')=1$ if $w(\Bg)=1$. Indeed, by the Argument Principle the curve
$\Bphi|D'$ intersects once every leave of $\di \Q$. Hence it has winding number
1 about this tube.
As the map 
$$g^{l_\la}: \Q\cap \pi^{-1} D'\rightarrow \V_0\cap \pi^{-1} D'$$ is a fiber bundles
diffeomorphism, it preserves the winding number.

Let us summarize the above discussion: 

\begin{lem}\label{generalized renormalization} 
Let $\Bg: \cup \V_i\ra\U$ be a  generalized  quadratic-like family over $D$.
Assume  it is proper and has  winding number 1. 
Then its generalized renormalization $\Bg': \cup\V_j'\ra \U'$ over $D'$ is also proper
and has  winding number 1.  
\end{lem}

\subsection{Central paracascades}\label{central paracascades}
In this section we will describe inductively tilings of the parameter plane according
to the types of generalized renormalizations. The initial tiling is constructed above
 in \S\ref{initial family}. Let us assume that we have already constructed a tiling
 $\DD^l$ of level $l$. The piece of this tiling containing a point $\la$ is denoted 
 by $D\equiv \Delta^l(\la)$. This piece comes together with   a holomorphic family
 $\Bg: \cup \V_i\ra\U$ of generalized quadratic like maps over $D$.

We will now subdivide  $D$ according to the combinatorics
of the central cascades of maps $g_\la$ (see \cite{L2}, \S\S 3.1, 3.6).
To this end let us first stratify the parameter values according to the length of their
central cascade. This yields a nest of parapuzzle pieces
$$D\supset D^{'}\supset\dots\supset D^{(N)}\supset\dots$$ 
For $\la\in D^{(N)}$, the map $g_\la$ has a central cascade 
\begin{equation}\label{central cascade}
V_\la^{(0)}\equiv  U_\la\supset   V_\la \equiv   
  V^{(1)}_\la\supset\dots\supset V^{(N)}_\la
\end{equation}
of length $N$, so that $g_\la 0\in V^{(N-1)}_\la\sm V^{(N)}_\la$.
Note that the puzzle pieces $V^{(k)}_\la$ are organized into the tubes
$\V^{(k)}$ over $D^{(k-1)}$ (with the convention that $D^{(-1)}\equiv D$). 

The intersection of these puzzle pieces, $\cap D^{(N)}$, is the little Mandelbrot set $M(\Bg)$
``centered" at the superattracting parameter value $c=c(\Bg)$ such that  $g_c(0)=0.$
Let us call $c$ the {\it center} of $D$. 

For a tube $\X$ over $D$ and a domain $\Delta\subset D$,
 let $\X|\Delta$ stand for the $\X\cap \pi_1^{-1}\Delta$. 
There is a Bernoulli map  
\begin{equation}\label{Bernoulli map}
\BG: \cup \W_{j}\ra \U
\end{equation}
 associated with the 
cascade \ref{central cascade} (see \cite{L3}, \S 3.6).  Here 
 the tubes  $\W_{j}$ over $D^{(N-1)}$
are the pull-backs of the tubes  $\V_i|D^{(k)}$, $i\not=0$,
 by the  covering maps 
\begin{equation}\label{coverings}
\Bg^k: (\V^{(k)} \sm \V^{(k+1)})|D^{(N-1)}\ra 
 ( \U\sm\V)|D^{(N-1)},\quad k=0,1\dots, N-1.   
\end{equation}

In the same way as in \S \ref{renormalization of families}, to  any string
$\j=(j_0,\dots,j_{l-1})$ corresponds the  tube over $D^{(N-1)}$,
 $$\W_\j=\{p\in\U|D^{(N-1)}: \BG^n p\in \W_{j_n},\; n=0,\dots, l-1\}. $$
Note that $\BG^l$ univalently maps each $ \W_\j$ onto $\U|D^{(N-1)}$. Thus 
$\W_\j$ contains a tube $\L_\j$ which is univalently mapped by $\BG^l$ onto
the central tube $\V^{(N)}$. These maps altogether form 
the first  landing map to $\V^{(N)}$ 
\begin{equation}\label{landing map}
\BT: \cup \L_\j\ra \V^{(N)}.
\end{equation}

For a $\la\in D^{(N-1)}\sm D^{(N)}$, let us now consider the itinerary $\j_\la$ of the
critical value $\phi(\la)\equiv g_\la(0)$
   through the tubes $W_j$ until its first return to $V^{(N)}$,
so that $\phi(\la)\in L_{\j_\la}$. Let
\begin{equation}\label{new tiles}
\Delta^{l+1}(\la)=\Bphi^{-1}\L_{\j_\la},\quad
\Q^{l+1}(\la)=\Bphi^{-1}\W_{\j_\la}.
\end{equation}

Thus  the annuli $D^{(N-1)}\sm D^{(N)}$ are tiled by  the parapuzzle pieces
$\Delta(\la)$ according
as the itinerary of the critical point through the Bernoulli scheme
(\ref{Bernoulli map})  until the first return to $V^{(N)}_\la$.
Altogether these tilings form the desired new subdivision  of $D$.
(Note however that the new tiles don't cover the whole domain  $D$:
the residual set consists of the Mandelbrot set $M(\Bg)$ and
of the  Misiurewicz parameter values for which the critical
orbit never returns back to $V_\la^{(N)}$, $\la\in D^{(N-1)}\sm D^{(N)}$. 

The affiliated holomorphic family over $\Delta(\la)$
 is defined as the first return map to $V^{(N)}_\la$. Its domain is obtained
by pulling back the tubes $\L_\j$ from (\ref{landing map})
by the double branched covering $\Bg: \V^{(N)}\ra \V^{(N-1)}|\Delta(\la)$,
and the return map itself is just $\BT\circ\Bg$. 

\subsection{Principal parapuzzle nest}\label{principal nest} 
Let us now summarize the above discussion. Given a quadratic-like family $\Bf$,
we consider the first tiling $\DD^0$ of a  Misiurewicz wake $O$ as described in 
\secref{first generalized family}. 
Each tile $\Delta\in \DD^0$ comes together with a generalized quadratic-like family
 $(\Bg_\Delta, \Bh_\Delta)$ over $\Delta$.

Now assume inductively that we have constructed the tiling $\DD^l$ of level $l$.
Then the tiling of the next level, $\DD^{l+1}$
 is obtained by partitioning each tile $\Delta\in \DD^l$
as described in \S\S\ref{renormalization of families},\ref{central paracascades}.

Let  $\Delta^l(\la)$ stand for the tile of $\DD^l$ containing $\la$, while
 $\Delta^l(\la)\subset Q^l(\la)\subset \Delta^{l-1}(\la)$
stand for another tile defined in (\ref{new tiles}).
Each tile $\Delta=\Delta^l(\la)$ contains a {\it central subtile}
 $\Pi^l(\la)=\Bphi^{-1}\V_0$ corresponding to the central return of the critical
point (here $\Bphi_\Delta=\Bg_\Delta(\bold 0)$).

 Let us then consider the sequence of  renormalized families $\Bg_{n}$
 over topological  discs
$\Delta^l(\la)$. We call the nest of topological discs  
 $\Delta^0(\la)\supset \Delta^1(\la)\supset \Delta^2(\la)\supset\dots$ 
{\it the principal parapuzzle nest}
of  $\la$. This nest is finite if and only if $\la$ is renormalizable.

Let $c_{n,\la}\in \Delta^l(\la)$ be the centers of the corresponding parapuzzle pieces.
We call them the  {\it principal superattracting approximations} to  $\la$.
If $\la$ is not renormalizable, then $c_{l,\la}\to \la$ as $l\to \infty$,
as $\diam \Delta^l(\la)\to 0$ (see the next section).   

The $\mod(\Delta^l(\la)\sm \Delta^{l+1}(\la))$
 are  called the {\it  principal parameter moduli} of $\la\in D$.

\section{Parapuzzle geometry}\label{parapuzzle geometry}

The following is the main geometric result of this paper:

\proclaim Theorem A. Let us consider a proper  DH quadratic-like family $\Bf$ 
with winding number 1 over $D$, and a Misiurewicz wake $O\subset D$.
Then for any $\la\in M(\Bf)\cap  O$, 
  $$\mod (\Delta^l(\la)\sm \Delta^{l+1}(\la))\geq Bl,\quad \text{and}\quad
 \mod(\Delta^l(\la)\sm \Pi^l(\la))\geq Bl,$$
where the constant  $B>0$ depends only on  $O$ and $\mod(\Bf)$. 

The rest of this section will be devoted to the proof of this theorem.

\subsection{A priori bound on parameter moduli }\label{a priori bounds}

  In this section we will show that the parameter annuli have definite moduli.
Given a holomorphic motion $h_\la$ and a holomorphic family of affine maps 
$g_\la: z\mapsto a_\la z+b_\la$, we can consider an \lq{affinely equivalent}\rq
  motion $g_\la\circ h_\la$.   In this way the motion can be normalized such that
any two points $z, \zeta\in U_*$ don't move (that is, $h_\la(z)\equiv z$ 
and $h_\la(\zeta)\equiv \zeta$
for $ \la\in D$). 
Let us start with a technical lemma:

\begin{lem}\label{round cylinders}
Let us have a holomorphic motion $h: (U_*,V_*)\ra (U_\la, V_\la)$ of a pair of nested
topological discs over a domain $D$. Assume that the maps 
$h_\la: (\di U_*,\di V_*)\ra (\di U_\la, \di V_\la)$ admit $K$-qc extensions 
$H_\la: (\C, U_*)\ra (\C, V_\la)$ 
(not necessarily holomorphic in $\la$ but with uniform dilatation $K$). Then there exists an
$M=M(K)$ such that if $\mod(U_*\sm V_*)>M$ then after appropriate normalization of the motion,
 there exists a round cylinder $D\times \{q<|z|<2 q\}$ embedded into
 $\U\sm\V$. 
\end{lem}   
  
\begin{pf} 
Let us normalize the motion in such a way that $0\in V_*$, and $h_\la(0)=0$.
Let $z_*$ be a point on $\di U_*$ closest to 0.
Normalize the motion in such a way that  $z_*=1$, and this point does not move either.
With this normalization, $V_*\subset D(0,\eps)$ where $\eps=\eps(m)\to 0$ as 
$m\equiv \mod(U_*\sm V_*)\to \infty$. 

Since the space of normalized $K$-qc maps is compact,  $|H_\la(\eps e^{i\theta}|)<\delta$,
where $\delta=\delta(\eps,K)\to 0$ as $\eps\to 0$, $K$ being fixed, and
$|H_\la (e^{i\theta})|>r$ where $r=r(K)>0$. It follows that the domain $\U$ contains
the round cylinder $D\times \{\delta<|z|<r\}$, and we are done. 
\end{pf}

\begin{cor}\label{parameter a priori estimate}
 Under the circumstances of \lemref{round cylinders},
let $\Bphi: D\ra\U$ be a proper
 analytic map with winding number 1. Let $D'=\Bphi^{-1} \V$. If $\mod(U_*\sm V_*)>M=M(K)$
  then $\mod(D\sm D')\geq\log 2$.
\end{cor}  

\begin{pf} By \lemref{round cylinders}, $\U\sm \V \supset D\times A$ where
$A= \{q<|z|<2q\}.$
  Let $Q=\Bphi^{-1} (D\times A)$. By the Argument Principle, 
$\phi=\pi_2\circ\Bphi$ univalently maps $Q$ onto $A$,
 so that $\mod(D\sm D')\geq \mod Q=\mod A=\log 2$.
\end{pf}

\subsection{A priori bound on dilatation}
Fix a point  $*$ in the Misiurewicz wake 
$ O$. Let us consider a
generalized quadratic-like  family $(\Bg: \cup \V_i\ra \U, \Bh)$
 over $D\equiv  \Delta^l(*)$. In what follows  we will use the notations of
 \secref{first generalized family}. Let $*\in D^{(N)}$. Let us lift the
holomorphic motion $\Bh$ from the cylinder  $\U\sm \V_0$ to the cylinders
$(\V^k\sm \V^{k+1})|D^{(N-1)}$ via the coverings (\ref{coverings}).  This provides
us with a holomorphic motion $\bold{H}$ of $\U\sm \V^{(N)}$ over $D^{(N-1)}$ with
the same dilatation as $\Bh$.

\begin{lem}\label{dilatation bound} 
The holomorphic motion  $\bold{H}$ of $\U\sm \V^{(N)}$ 
 over $\Delta^{l+1}(*)$ (see (\ref{new tiles}))
 has a uniformly bounded dilatation, depending only on the choice of the Misiurewicz  wake 
$O$ and $\mod(\Bf)$.
\end{lem}  

\begin{pf} Let us go back to (\ref{new tiles}).
 By \cite{L2}, $\mod (W_{j_*}\sm L_{j_*})$ is big for big  $l$.
By \cite{L3}, for $\la\in Q^{l+1}(*)$, 
there is a $K$-qc pseudo conjugacy
$\psi_\la: (U_*, \cup V_{i,*})\ra (U_\la, \cup V_{i,\la })$, with $K$ depending
only on the choice of wake $O$ and $\mod(\Bf)$.
 Hence Corollary \ref{parameter a priori estimate}
can be applied. We conclude that 
\begin{equation}\label{a priori bound}
\mod(Q^{l+1}(*)\sm\Delta^{l+1}(*))\geq \log 2
\end{equation}
for $l$ sufficiently big (depending on $O$ and $\mod(\Bf)$).

The lift of $\bold{H}$  to the  landing tubes (\ref{landing map}) yields a 
$K$-qc holomorphic motion
 of ($\U\sm \V^{(N)},  \cup \L_{\bar j})$ over $D^{(N-1)}$.
By the Extension Lemma, this motion can be extended through $\V^{(N)}$.
Let us  keep the same notation $\bold{H}$ for this motion. 

By the Quasi-Conformality Lemma and (\ref{a priori bound}),
  $\bold{H}$  is $K$-qc over $\Delta^{l+1}(*)$, 
with an absolute $K$ provided $l$ is big enough.
But the holomorphic motion $\Bh_{l+1}$ on  $\U^{l+1}\sm \V^{l+1}_{0}$
is the lift of $\bold H$ on $\V^{(N-1)}\sm \L_{\bar j_*}$ over $\Delta^{l+1}(*)$
 via the fiberwise  analytic  double covering
 $$\Bg_l: \U^{l+1}\sm \V^{l+1}_0\ra \V^{(N-1)}\sm \L_{\bar j_*}. $$
 Hence $\Bh_{l+1}$  is also $K$-qc.
\end{pf} 

\subsection{Proof of Theorem A } We are now prepared 
to complete the proof:
$$\mod (\Delta^l\sm \Delta^{l+1})
\geq \mod(W_{\bar i_*}\sm L_{\i_*})\geq Bl.$$
The first  estimate in the above row  follows from Lemma \ref{dilatation bound}
and the Remark following  Corollary \ref{transversal qc structure}.
 The last estimate is the main result of
\cite{L2}.

For  the same reason,
$$\mod(\Delta^l\sm \Pi^l)\asymp \mod(U^l_*\sm V^l_{0,*})\geq Bl.$$

\section{Application to the measure problem}\label{measure}

In this section we will apply the previous results to the real 
quadratic family $P_c: z\mapsto z^2+c$, $c\in \R$. Let $\NN\RR$ stand for the set of
non-renormalizable real parameter
values $c\in [-2,-3/4)$. Note that all periodic points of
 the $P_c: z\mapsto z^2+c$, $c\in\NN\RR$, are repelling. 
Indeed, the interval $[-3/4, 1/4]$ where
$P_c$ has a non-repelling fixed point is excluded, while maps with
non-repelling cycles of higher period are renormalizable.

Let $\NN\CC$ stand for the set of parameter values $c\in\NN\RR$ 
such that the principal nest of $P_c$ constants only finitely
many non-trivial (i.e., of length $>1$) central cascades. 

\begin{thm}\label{measure zero}
\begin{itemize}
\item The set $\NN\RR$ has positive measure;
\item The set $\NN\CC$ has full  Lebesgue measure in $\NN\RR$.
\end{itemize}
\end{thm}

\noi{\bf Remarks.} 1.
The former (positive measure) result is known (see \cite{BC}, \cite{J}).
The latter (full measure) is new.

2. The corresponding statements concerning at most finitely renormalizable
  parameter values are derived from the above statements by considering quadratic-like
families associated with little copies of the Mandelbrot set. 

3. By the result of  Martens \&  Nowicki \cite{MN}  together with \cite{L5},
$P_c$ has an absolutely continuous invariant measure for any $c\in \NN\CC$. 
Altogether these yield Theorem B stated in the Introduction. \smallskip

\noindent{\it Proof of Theorem \ref{measure zero}.}
Let $d$ stand for the real tip of the little Mandelbrot set
attached to the main cardioid (i.e.  $P_d^3(0)=\alpha$).
As  all parameter values $c\in [d, -3/4)$ are renormalizable,
we can restrict ourselves to the interval $[-2,d)\supset \NN\RR$. This interval
belongs to the Misiurewicz wake $O$ attached to $d$.

Given measurable sets $X,Y\subset \R$, with $\length(Y)>0$, let
 $\dens(X|Y)$ stand for the $\length(X\cap Y)/\length(Y)$.

We will now restrict all tilings $\DD^l$ constructed above to the real line,
without change of notations.
We will use the same notation, $\DD^l$, for the union of all pieces of $\DD^l$.
 For every
 $\Delta=\Delta^{l}(\la)\in \DD^{l}$,
let us consider the  central piece  $\Pi\subset \Delta$ 
 corresponding to the central  return of the critical point. 
By Theorem A, $\dens(\Pi|\Delta) \leq Cq^l$  
 for absolute  $C>0$ and $q<1$. Let $\Gamma^l$
 be the union of these central pieces.
Summing up over all $\Delta\in \DD^l$, we conclude that 
\begin{equation}\label{exp decay}
\length(\Gamma^l)\leq \dens(\Gamma^l|\DD^{l})\leq Cq^l
\end{equation}
(the whole interval is normalized so that its length is equal to 1).

It follows that for $l$ sufficiently big,
$$\dens(\bigcup_{k\geq 0}\Gamma^{l+k}|\DD^l)\leq C_1 q^n<1,$$
which means that with positive probability central returns will never occur again. 
This proves the first statement.

To prove the second one just notice that (\ref{exp decay}) together with the Borel-Cantelli
Lemma yield that infinite number of central returns occurs with zero probability.
\QED

\end{document}